\documentclass{article} \usepackage{amssymb}  
 \newcommand{\La}{\Lambda}
\newcommand{\R}{{\mathbb R}} \newcommand{\Z}{{\mathbb Z}} \newcommand{\N}{{\mathbb N}}
 \newcommand{\C}{{\mathbb C}}  

  \newcommand{\cc}{\xi}\newcommand{\K}{\mathcal{K}} \newcommand{\Ss}{\mathcal{S}} \newcommand{\Q}{\mathcal{Q}} \newcommand{\B}{\mathcal{B}}

\begin{document}
\begin{large}

\title{On multi--dimensional sampling and interpolation\thanks{Published in Anal. Math. Phys. 2 (2012), no. 2, 149--170}}

\author{Alexander Olevskii\thanks{The first author is supported in part by the Israel Science
        Foundation} \ and Alexander Ulanovskii\thanks{The second author  is supported by an ESF-HCAA grant}}

\date{} \maketitle

\noindent A.O.: School of Mathematics, Tel Aviv University\\ Ramat Aviv,  69978 Israel\\ E-mail:
 olevskii@post.tau.ac.il


\noindent A.U.: Stavanger University,  4036 Stavanger, Norway\\ E-mail: Alexander.Ulanovskii@uis.no

 \begin{abstract}

    The paper discusses sharp sufficient conditions for interpolation and sampling for functions of $n$ variables with convex spectrum.   When $n=1$, the classical theorems of Ingham and Beurling state  that the critical values in the estimates from above     (from below) for the distances between interpolation (sampling) nodes are the same.
    This is no longer true for  $n>1$.    While the critical value for sampling sets remains constant,
    the one for interpolation grows linearly  with the dimension.
\end{abstract}

  \section{Introduction}

      The paper  discusses sampling and interpolation problems in $\R^n$.
      We focus on sharp sufficient conditions  for the Paley--Wiener and Bernstein spaces
      of entire functions with convex spectrum.

    \subsection{Paley--Wiener and Bernstein Spaces}

     {\bf Definition 1.1} {\it Let $\Ss\subset \R^n$ be a bounded set of  positive measure.
   The Paley--Wiener space $PW_\Ss$ consists of all Fourier transforms
   $$
   \hat F(t)=\frac{1}{(2\pi)^{n/2}}\int_{\R^n}e^{-i t\cdot x}F(x)\,dx
   $$
   of
   functions $F\in L^2(\R^n)$ which vanish a.e. outside $\Ss$.}

   \medskip
  Equipped with the $L^2$--norm, $PW_\Ss$ is a Hilbert space.

 \medskip\noindent{\bf Definition 1.2} {\it Let $\Ss\subset \R^n$ be a compact set.
 The  Bernstein space  $ B_\Ss$ consists of all continuous
 bounded functions  on $\R^n$ which are the Fourier transforms of distributions supported by $\Ss$.}

 \medskip

  Equipped with the $L^\infty-$norm, $B_\Ss$ is a Banach space.

Let $\Ss$ be a closed convex symmetric body in $\R^n$.
      We denote by $\Vert x \Vert_\Ss=\min\{r\geq 0: x\in r\Ss\}$ the norm generated by  $\Ss$, and by
      $\Ss^o:=\{x\in \R^n: x\cdot y\leq 1, y\in\Ss\}$  the polar body of $\Ss$.

  When $\Ss$ is a closed convex symmetric body, the Paley--Wiener and Bernstein spaces admit  a simple characterization:

\medskip\noindent{\bf Theorem A} {\sl (i)  $PW_\Ss$ consists of all entire functions $f$ which belong to $L^2(\R^n)$ and satisfy
\begin{equation}\label{pw}
|f(u+iv)| \leq C e^{\Vert v\Vert_{\Ss^o}}, \ u,v\in\R^n,
\end{equation}with some constant  $C$.

(ii)  $B_\Ss$ consists of all entire functions $f$ which satisfy (\ref{pw}) with some constant $C$. }

\medskip

The first part is a multi-dimensional version of the classical Paley--Wiener theorem, see for example \cite{sw}, ch. 4.

     The second part can be deduced easily from the first one  by  convolving the distribution $\hat f$ with a compact
     $C^\infty$-approximative unity.

Observe that  in (\ref{pw}) one can take $C=\sup_{t\in \R^n}|f(t)|$.

\subsection{Sampling and Interpolation}


     {\bf  Definition 1.3}. (i) {\sl  A set $\La\subset\R^n$ is called uniformly discrete if}
\begin{equation}\label{ud}
              \delta(\La):=\inf_{\lambda,\lambda'\in\La,\lambda\ne\lambda'}|\lambda-\lambda'|>0.
\end{equation}

  (ii) {\sl A set $\La\subset\R^n$ is called relatively dense  if}
\begin{equation}\label{rd}
              \rho(\La):=2\sup_{x\in\R^n}\inf_{\lambda\in\La}|x-\lambda|<\infty.
\end{equation}

\medskip


 Clearly, $\rho(\La)=\delta(\La)$ when $\La$ is an arithmetic progression in $\R$.


       The {\it Sampling Problem} deals with possibility of reconstruction of a continuous-time signal $f$
       with spectrum in $\Ss$  from its samples    (values of $f$ on $\La$).             We shall be interested  in the stable   sampling.

         \medskip\noindent {\bf  Definition 1.4}.  {\sl     (i)
     $\La$ is called a stable sampling set (SS) for $PW_\Ss$ if
\begin{equation}\label{fp}
               \Vert f\Vert^2_2 \leq  C \Vert f|_\La\Vert^2_2 , \ f\in PW_\Ss,
\end{equation}     where the constant $C$ does not depend on $f$.

(ii) $\La$ is called an SS for $B_\Ss$ if \begin{equation}\label{sb}
\Vert f\Vert_\infty\leq C\Vert f|_\La\Vert_\infty, \ f\in B_\Ss,
\end{equation}where the constant $C$ does not depend on $f$. }

\medskip
Here we used the notation
$$
\Vert f\Vert^2_2:=\int_\R|f(t)|^2\,dt, \Vert f|_\La\Vert^2_2:=\sum_{\lambda\in\La} |f(\lambda)|^2
$$
and
$$
\Vert f\Vert_\infty:=\sup_{t\in\R}|f(t)|, \ \Vert f|_\La\Vert_\infty=\sup_{\lambda\in\La}|f(\lambda)|.
$$

It is well-known that  a sampling set is always relatively dense.

When $\La$ is uniformly discrete, the inverse inequality in (\ref{fp}) holds  true (see \cite{ya}, Theorem 2.17).  This implies that  $\La$ is an SS  for $PW_\Ss$ if and only if  the exponential system \begin{equation}\label{ex}E(\La):=\{e^{i\lambda\cdot x}{\bf 1}_\Ss(x),\lambda\in\La\}\end{equation} is a frame in $L^2(\Ss)$.

  The {\it Interpolation Problem} is in a way dual to the
         sampling one.

\medskip
\noindent {\bf  Definition 1.5}. {\sl
 (i) A countable set $\La$ is called an  interpolation set (IS) for $PW_\Ss$  if for every $c(\lambda)\in l^2(\La)$ there exists $f\in PW_\Ss$ satisfying $f(\lambda)=c(\lambda),\lambda\in\La$.

   (ii) A countable set $\La$ is called an  IS    for $B_\Ss$ if for every $c(\lambda)\in l^\infty(\La)$ there exists $f\in B_\Ss$ satisfying $f(\lambda)=c(\lambda),\lambda\in\La$.}

\medskip

It is well-known that  an interpolation set is always uniformly discrete.

         In geometric language, $\La$ is an IS for $PW_\Ss$ if and only if
          the exponential system $E(\La)$ is a Riesz basis in its
          linear span,           that is the two-side estimate holds:
\begin{equation}\label{ip}
  C_1\Vert c\Vert_2\leq \Vert \sum_{\lambda\in\La} c(\lambda)e^{i\lambda\cdot x}\Vert_{ L^2(\Ss)} \leq C_2\Vert c\Vert_2,
\end{equation} for every finite sequence  $c=c(\lambda),$ where $C_j=C_j(\Ss,\La), j=1,2,$ are positive constants.

The right hand-side inequality in (\ref{ip}) is true whenever $\La$ is a uniformly discrete  set and $\Ss$ is a bounded set (\cite{i}, Theorem 2, \cite{ya}, p.135).

Notice that if $\La$ is an IS for $PW_\Ss$, the solution to the interpolation problem $f(\lambda)=c(\lambda), c\in l^2(\La),$ can be chosen with the additional
   requirement
   \begin{equation}\label{si}
                 \Vert f\Vert_{2}\leq  C \Vert c\Vert_2,
\end{equation}where $C>0$ is independent on $c$. A similar property is valid for the interpolation in $B_\Ss.$

    \subsection{ Classical results}

  It is a  classical problem  to describe sampling and interpolation sets for a given  space.

         Let $\Ss$ be an interval in $\R$ and let  $|\Ss|$ denote its length.

     \medskip\noindent{\bf Theorem B}   {\sl
        If   $\rho(\La) < \frac{2\pi}{|\Ss|}$ then $\La$ is an SS for $PW_\Ss.$}

        \medskip\noindent{\bf Theorem C}   {\sl
        If  $\delta(\La) > \frac{2\pi}{|\Ss|}$   then $\La$ is an IS for $PW_\Ss$.}

        \medskip

         Here $\rho(\La)$ and $\delta(\La)$ are defined in (\ref{rd}) and (\ref{ud}), respectively.

         These results hold also for the Bernstein space $B_\Ss$.

        Theorem B is a one-dimensional version of Beurling's result from \cite{b3}. Theorem C is due to Ingham \cite{i}.

         More general conditions for sampling and interpolation on a single interval were given by Beurling \cite{b1}, \cite{b2}
        and Kahane \cite{k2} in terms of appropriate uniform densities (the lower
        and upper uniform densities). These results are based    on the entire functions theory.
         Landau \cite{la}  extended the {\it necessary} density conditions in these  results to general bounded spectra in $\R^n$.
                   However, the {\it sufficient} density conditions fail
                   for disconnected spectra in $\R$.  This is the place where the arithmetics of $\La$ comes into
         the play.

   In $\R^n, n>1,$ the situation is even "worse".

    \subsection{ Several variables}

          There is a fundamental difference between sampling and           interpolation of functions of one and several variables:
           In dimension $n = 1$ the zeros of an entire function $f$ are
          discrete, and there is a precise connection between the
          asymptotic behavior/density of the zeros and the growth of $f$.
           In several dimensions, the zero sets are analytic manifolds.
           That is why sufficient conditions for sampling and
          interpolation cannot be given in terms of the uniform densities. See also discussion in Seip \cite{seip}, p. 122.

            Theorems B and C turns out to be more effective from           this point of view.

          The present paper  studies {\it sufficient} conditions for sampling and interpolation in the
          situation when the spectrum is a convex body in $\R^n$.

          The paper is organized as follows.           In sec.2 we consider connection between
          sampling/interpolation problems in Paley--Wiener and           Bernstein spaces.

         In sec.3 we discuss the sampling/interpolation on lattices.          We also  describe small perturbations of
         lattices which provide universal  SS and IS of optimal density.
         This construction in the one-dimensional case was presented
         in \cite{ou0} and \cite{ou}.

Sec. 4 is devoted to sampling.           We suggest an alternative approach to Beurling's  Theorem B in $\R^n$. This approach allows us to avoid
         a non-trivial balayage technique  introduced by Beurling and present the result
         in a slightly stronger form.          It also makes quite visible its one-dimensional nature.

           The latter is not the case for the interpolation problem.
          In sec.5 we consider this problem in a general context,    when the spectrum is a ball in some metric
          in $\R^n$,  and the distances between interpolation nodes     are measured in another metric.
          We give a necessary condition for interpolation.      The proof is based on Minkovski's lattices.
          Our sufficient condition involves   concentration property of functions and their
          Fourier transform.     We illustrate the results in the case of $l^p$--metrics,
          when these two conditions provide estimates which are    asymptotically sharp with respect to the dimension.
           In particular, this shows a contrast in the behavior    of the critical Beurling and Ingham bounds  for large dimension.

In what follows we  denote by 
 $|\Ss|$ the ($n-$dimensional) measure of a set $\Ss\subset\R^n$,  $|x|$ the Eucledian norm of $ x\in \R^n$, $\B:=\{x\in\R^n:|x|\leq 1\}$ the closed unit ball in $\R^n$ and  $$r\Ss=\{rx: x\in \Ss\},  \Ss+\K=\{x+y: x\in\K, y\in \Ss\}, r>0,\K\subset\R^n.$$

  \section{Connection between sampling/interpolation in PW- and B-spaces}

\noindent{\bf Theorem 2.1} {\sl  Suppose $\Ss$ is a compact, $\La$ is a uniformly discrete set in $\R^n$ and $\epsilon>0$.


(i) If $\La$ is an SS for $PW_{\Ss+\epsilon\B}$ then it is a SS for $B_{\Ss}$.

(ii) If $\La$ is an SS for $B_{\Ss+\epsilon\B}$ then it is a SS  for $PW_{\Ss}$.

(iii)  If $\La$ is an IS for $PW_{\Ss}$ then it is an IS for $B_{\Ss+\epsilon\B}$.

 (iv) If $\La$ is an IS for $B_\Ss$ then it is an IS for $PW_{\Ss+\epsilon\B}$}

\medskip
 This result is basically known but we have not found it in this form in the literature.

           The proofs of (i), (iii) and (iv) are easy.
           To proof  (ii)
 one may use Beurling's Linear Balayage Operator (see \cite{b2}, pp.348--350,  \cite{b3}, pp.306--310). We  present a new elementary proof of (ii).

 \medskip\noindent{\bf Proof}. (i)   Suppose $\La$ is not an SS for $B_{\Ss}$. This means that there are functions $g_j\in B_\Ss$ satisfying $|g_j(x_j)|=1,$ for some $x_j\in\R^n$ and $\Vert g_j|_\La\Vert_\infty<1/j.$ Fix any function $\Phi\in PW_{\epsilon\B}$ satisfying $\Vert \Phi\Vert_\infty=\Phi(0)=1$. Then
the function $f_j(x):=\Phi(x-x_j) g_j(x)$ belongs to $PW_{\Ss+\epsilon\B}$ and satisfies
$ |f_j(x_j)|=1$. From this, by Bernstein's inequality for entire functions of exponential type, we get $\Vert f_j\Vert_2\geq K$, where $K>0$ depends only on the diameter of $\Ss$. On the other hand, by the inverse inequality of (\ref{fp}), we have
$$
\Vert f_j|_\La\Vert_2\leq \frac{1}{j}\Vert \Phi(x-x_j)|_{\La}\Vert_2\leq \frac{1}{jC}\Vert \Phi\Vert_2\to 0,\ j\to\infty.
$$
This  contradicts to  (\ref{fp}), and so $\La$ is an SS for $B_{\Ss}.$

(ii)  We may assume that $\Ss\subseteq(-\pi,\pi)^n$. Take a function $\Phi\in PW_{\epsilon \B}$ with $\Phi(0)=1$ and set
$$
K:=\sup_{x\in [0,1]^n}\sum_{m\in\Z^n}|\Phi(m-x)|^2<\infty.
$$
For every $f\in PW_{\Ss}$ we have $|f(m)|\leq \sup_{x\in\R^n}|f(x)\Phi(m-x)|$ and $f(x)\Phi(m-x)\in B_{\Ss+\epsilon\B}, \ m\in\Z^n.$  Using (\ref{sb}), we obtain
$$
\Vert f \Vert_{2}^2=\sum_{m\in\Z^n}|f(m)|^2\leq \sum_{m\in\Z^n} \sup_{x\in\R^n}|f(x)\Phi(m-x)|^2\leq
$$$$ C^2\sum_{m\in\Z^n} \max_{\lambda\in\La}|f(\lambda) \Phi(m-\lambda)|^2\leq $$$$C^2\sum_{m\in\Z^n}\sum_{\lambda\in\La} |f(\lambda)\Phi(m-\lambda)|^2
\leq C^2K\Vert f|_\La\Vert_2^2.
$$

(iii)
Take a function $\Phi\in PW_{\epsilon\B}$ satisfying $\Phi(0)=1$ and \begin{equation}\label{v}|\Phi(x)|\leq \frac{K_1}{(1+|x|)^{2n}}, \ x\in \R^n.\end{equation} Choose functions $\varphi_\lambda\in PW_\Ss, \lambda\in\La,$  satisfying $\sup_{\lambda\in\La}\Vert \varphi_\lambda\Vert_2<\infty$,
$\varphi_\lambda(\lambda)=1$ and $\varphi_\lambda(\lambda')=0, \lambda\in\La, \lambda'\ne\lambda.$ Clearly, we have
$\Phi(x-\lambda)\varphi_\lambda(x)\in PW_{\Ss+\epsilon\B}$.

Given any sequence $\cc(\lambda)\in l^\infty(\La)$, we wish to find a function $f\in B_{\Ss+\epsilon\B}$ satisfying $f(\lambda)=\cc(\lambda), \lambda\in\La$. Observe that  $$\Vert \varphi_\lambda\Vert_\infty^2\leq \frac{|\Ss|}{(2\pi)^n}\Vert \varphi_\lambda\Vert_2^2<\infty.$$ This, (\ref{v}) and the uniformly discreteness of $\La$ shows that the function
$$
f(x):=\sum_{\lambda\in\La}\cc(\lambda)\Phi(x-\lambda)\varphi_\lambda(x)
$$ belongs to $B_{\Ss+\epsilon\B}$. Clearly, this function
solves the interpolation problem  above.

 The proof of (iv) is similar to (iii).
         Take a function $\Phi$ satisfying (\ref{v}) and functions  $\varphi_\lambda\in B_\Ss, \lambda\in\La,$ satisfying $\varphi_\lambda(\lambda)=1$, $\varphi_\lambda(\lambda')=0, \lambda'\ne\lambda$, and set $K:=\sup_\lambda\Vert\varphi_\lambda\Vert_\infty<\infty.$
        It suffices to verify that the function $f$ defined above belongs to $L^2(\R^n)$. Indeed, from the right hand-side inequality in (\ref{ip}), we obtain:
 $$\Vert f\Vert_2^2\leq K\Vert \sum_{\lambda\in\La}|\cc(\lambda) \Phi(x-\lambda)|\, \Vert_2^2=K\int_{\epsilon\B}| \sum_{\lambda\in\La}|\cc(\lambda)|e^{i\lambda\cdot t}\hat\Phi(t)
|^2\,dt
$$
$$
\leq K \Vert\hat\Phi\Vert_\infty^2\Vert \sum_{\lambda\in\La}|\cc(\lambda)|e^{i\lambda\cdot t} \Vert_{L^2(\epsilon\B)}^2<\infty,
$$   which finished the proof.

\section{Lattices and their perturbations}

\subsection{Lattices}


Let us start with the simple lattice $\Z^n$.
  The trigonometric system $E(\Z^n)=\{\exp(ik\cdot t), k\in\Z^n\}$
  forms an orthogonal basis in $L^2$ on the torus $[0,2\pi]^n$.
     Using periodicity argument one can obtain

\medskip\noindent
{\bf Proposition 3.1}  {\sl
(i)  $\Z^n$ is an SS for $PW_\Ss$ if and only if
     the translates $\Ss + 2\pi m , m\in \Z^n,$ satisfy the
    "packing" property: $$|\Ss\cap(\Ss+ 2\pi m)|=0, \ m\in \Z^n, m\ne 0.$$

    (ii) $\Z^n$ is an IS for $PW_\Ss$ if this set of translates satisfies
     the "covering"  property: $$\Ss+ 2\pi\Z^n=\R^n.$$ }

    Now let $\La$ be a general lattice:  $\La=T\Z^n$,
 where $T:\R^n\to\R^n$ is an invertible linear operator.   Denote by $\La^*:=(T^*)^{-1}\Z^n$ the dual lattice and by
det$(\La)=|T ([0,1]^n)|$ the determinant of $T$.

    The following proposition is straightforward:

   \medskip\noindent
{\bf Proposition 3.2}  {\sl $\La$ is an SS  (IS) for $PW_\Ss$ if and only if
             $\Z^n$ is SS (IS) for $PW_{T^*(\Ss)}$}.

\medskip
    The two propositions above imply

         \medskip\noindent
{\bf Proposition 3.3}  {\sl Let $\La=T\Z^n$ and $\Ss\subset\R^n, n\geq1$.

(i) $\La$ is an SS for $PW_\Ss$ if and only if it satisfies the packing property
$$
|\Ss\cap(\Ss+ 2\pi \lambda^*)|=0, \ \ \lambda^*\in\La^*, \ \lambda^*\ne0.
$$

 (ii) $\La$ is an IS for $PW_\Ss$ if and only if it satisfies the covering property}
 $$
 \Ss+2\pi\La^*=\R^n.
 $$

 \medskip

Below we use also the concept of the uniform density:

\medskip\noindent{\bf Definition 3.1} {\sl
We say that a set $\La\subset\R^n$  possesses a uniform density $D(\La)$ if}
$$
\frac{\#(\La\cap (x+r\B))}{|r\B|}= D(\La)(1+o(1)), \ r\to\infty \mbox{ uniformly on } x.
$$
Clearly, for every lattice $\La$ we have $D(\La)=1/$det$(\La)$.

We remarked in sec. 1.4 that in several dimensions neither sampling nor interpolation property of a set $\La$ can be formulated in terms of its uniform density. This phenomenon can be observed already for lattices. The following  is an easy consequence of Proposition 3.3:

\medskip\noindent
{\bf Corollary 3.1}. {\sl Let $n\geq 2$ and $\Ss\subset\R^n $ be a bounded set of positive measure.
For every $\epsilon>0$ there is a lattice $\Lambda_1\subset\R^n$ with $D(\La_1)>1/\epsilon$ which is not a SS  for  $PW_\Ss$, and  a lattice $\Lambda_2\subset\R^n$ with $D(\La_2)<\epsilon$ which is not an IS for  $PW_\Ss.$ }

\medskip

By Theorem 2.1, the same result holds for Bernstein spaces $B_\Ss$.

Observe that a similar phenomena occurs for the uniqueness sets. A set $\La$ is called a set of uniqueness for $PW_\Ss$ if there is no non-trivial function $f\in PW_\Ss$ which vanishes on $\La$. When $\Ss$ is an interval in $\R$, it is well--known that every set $\La\subset\R$ for which $D(\La)$ exists and satisfies $D(\La)<|\Ss|/2\pi$, is not a set of uniqueness for $PW_\Ss$. Such a result is no longer true in several dimensions. The following is an easy consequence of \cite{u}, Corollary 1:

\medskip\noindent
{\bf Corollary 3.2}. {\sl Let $n\geq 2$ and $\Ss\subset\R^n $ be a bounded set of positive measure.
For every $\epsilon>0$ there is a set $\Lambda\subset\R^n$ with $D(\La)<\epsilon$ which is  a uniqueness set for  $PW_\Ss$. }
\medskip

Observe that the set $\La$ in Corollary 3.2 can be chosen as an arbitrarily small perturbation of the lattice $\La_2$ from Corollary 3.1.

\subsection{Perturbed lattices  and Riesz bases}

 Recall that a family of vectors $\{f_k\}$ in a Hilbert space $H$       is called a Riesz basis if it is obtained from an orthonormal
      basis by the action of a linear isomorphism of the space.        An exponential system $E(\La)$ defined in (\ref{ex}) is a Riesz basis in $L^2(\Ss)$
      if and only if the set $\La$ is both SS and IS for $PW_\Ss$.        It is well--known that if $E(\La)$ is a Riesz basis in $L^2(\Ss)$
       then $D(\La)= |\Ss|/(2\pi)^n$ (see \cite{la}).

                    Observe that   the problem of existence of exponential Riesz bases in $L^2$ on an interval in $\R$ is well understood (see \cite{seip}). However,  except for the product domains like cube, there are very few examples of {\sl convex sets} $\Ss$ in $\R^n, n>1,$ for which it is known that a Riesz basis of exponentials in $L^2(\Ss)$ exists. For example, it is an open problem if such a basis exists in $L^2(\B)$, where $\B$ is the unite ball. We mention \cite{lr}, where exponential Riesz basis are constructed for convex polygons in $\R^2$, symmetric with respect to the origin.

      Given countable sets  $\La=\{\lambda(k)\}$ and $\La'=\{\lambda'(k)\}$ in $\R^n$, we say that $\La'$ is
      an $\epsilon-$perturbation of $\La$ if $\sup_k |\lambda(k)-\lambda'(k)| < \epsilon$
      for some enumerations of the sets.

      Small perturbation of a lattice  may improve essentially its sampling/interpolation properties.
      In \cite{ou0} and \cite{ou}  we proved that certain arbitrary small perturbations
      $\La$ of the lattice $\Z$  produce exponential systems $E(\La)$
      which constitute Riesz basis on every set $\Ss ,|\Ss|=2\pi$, which is a
      finite union of intervals whose endpoints belong to $\pi\mathbb Q$.
      We call such systems "universal"  Riesz bases.        It is mentioned in \cite{ou} that a similar result holds in $\R^n, n>1$.
      Here we give more details for this case.


By Proposition 3.2, we may restrict ourselves to the perturbations of $\Z^n$.
Denote by $\Omega$ the collection of all sets $\Ss$ which admit a representation

$$
\Ss=x+2^{-k}\bigcup_{j=1}^{2^{kn}}([0,2\pi]^n+2\pi m_j),
$$
for some $x\in\R^n$, $k\in\N$ and different vectors $m_1,...,m_{2^{kn}}\in\Z^n$. Observe that $|\Ss|=(2\pi)^n$ for every $\Ss\in\Omega$.

\medskip
\noindent {\bf Theorem 3.1} {\sl For every  $\epsilon>0$, there is an $\epsilon-$perturbation $\La$ of
 $\Z^n$ such that  the system
 $E(\Lambda)$ is a Riesz basis in $L^2(\Ss)$ for every
$\Ss\in\Omega$.}

\medskip

One can construct universal Riesz basis on a wider than $\Omega$ collections of sets as well.

Observe that
 $\Omega$ is sufficiently large in the sense that that every compact set $\Ss'$ of measure less than $(2\pi)^n$ can be covered
by a set $\Ss\in\Omega$. Also, every open set $\Ss'$ of measure $>(2\pi)^n$ contains a  set $\Ss\in\Omega$. Hence, By Theorem 3.1 and Proposition 3.2 we get

 \medskip\noindent{\bf   Theorem  3.2} {\sl For every lattice $\Gamma\subset\R^n$ and every $\epsilon>0$ there is an $\epsilon-$perturbation $\La$ of $\Gamma$ such that

(i) $\La$ is an SS for $PW_\Ss$ for each compact set $\Ss$, $|\Ss|<D(\La)(2\pi)^n$;

(ii)  $\La$ is an IS for $PW_\Ss$ for each open set  $\Ss$, $|\Ss|>D(\La)(2\pi)^n$.}

\medskip

         In  1D-case this was proved in \cite{ou0}, \cite{ou}.       We call such sets $\La$ "universal" sampling /interpolation
       sets.           For another construction of universal SS and IS,        based on Meyer's quasi-crystals, see \cite{m1}-\cite{m3}.

\subsection{Lemmas}

Set$$  I := 2\pi [0,1]^n,\   I(m):=I + 2\pi m, \ m\in \Z^n.$$

      \medskip\noindent{\bf Lemma 3.1.} {\sl
       Let $\Ss := \cup_{l=1}^k I(m_j)$, where $m_j$ are $k$ different
       vectors in $\Z^n$ and let
$           \La:= \cup_{l=1}^k(\Z^n + u_l),   $ where $u_l$ are   $k$  vectors in $\R^n$.
       Consider $k \times k$ matrix  $A:= (\exp(2\pi u_l\cdot m_j))$.
       Then $E(\La)$ is a Riesz basis in $L^2(\Ss)$ if and only if}
\begin{equation}\label{x}
                     \mbox{det}\, A \neq 0.
\end{equation}

       For $n=1$ this lemma is well-known, see   \cite{BK},  \cite{H}, p. 143--145.
       The main idea appeared already in \cite{L1}.

\medskip\noindent{\bf        Proof.}
        Consider the operator  $R$,
        $$
 R c=\sum_{\lambda\in \La}c(\lambda)e^{i\lambda\cdot t}{\bf 1}_\Ss(t), \ c\in l^2(\La).
$$
        Clearly it is  bounded  from $l^2(\La)$ to $L^2(\Ss)$.
        We need to check that this operator is one-to-one.
         Take any function $g\in L^2(\Ss)$ and write it as a sum
        of translates:
$$
                  g(t):= \sum_{j=1}^k g_j (t- 2\pi m_j) ,
$$
         where $g_j$ are functions supported by $I$.
         We have also: $$c=\bigcup_{l=1}^k c(l), \ c(l) \in l^2(\Z^n+u_l) ,$$
         so that $$                 Rc = \sum_{l=1}^k f_l  ,$$ where
$$
                 f_l(t)= \sum_{p\in\Z^n} c_l(p)e^{i(p+u_l)\cdot t}{\bf 1}_\Ss(t).
$$
          For $t=2\pi m_j + x , x\in I(0),$ we have:
$$
                f_l(t)= f_l(x) e^{2\pi iu_l \cdot m_j}.
$$
         It follows that the equation $Rc=g$ can be
         re-written as the system of linear equations
         with respect to ${f_l }$ with matrix $A.$
         Condition (\ref{x}) implies uniqueness of the
         solution in $[L^2(I)]^k.$
         Expanding each $f_l$ in the Fourier series
         in the orthogonal basis $E(\Z^n + u_l)$,
         we get a unique vector $c\in l^2(\La)$
         which solves the equation $Rc=g.$
         The necessity of (\ref{x}) is also clear.

\medskip\noindent{\bf  Lemma 3.2} {\sl Let $v_j,  j=1,2,...,k,$ be         different vectors in $\C^n.$
        For $u_l \in \C^n, l=1,2,...,k$, denote
$$
              h(u_1,u_2, ... u_k):= \mbox{det} (e^{ 2\pi u_l\cdot v_j}).
$$
         Then $h$ is a non-trivial entire function in $\C^{kn}$.}

\medskip
          Indeed, it is clear that the determinant above is an entire function of $u_l$ and that the $u_l$  may be chosen so that it becomes a Vandermonde determinant different from zero.

\medskip

\noindent{\bf Lemma  3.3} (A version of Paley--Wiener Stability Theorem) {\sl
Assume $E(\Gamma)$ is a Riesz basis in $L^2(\Ss)$ for some compact  $\Ss\subset\R^n$.  Then there exists $\delta>0$ such that for every $\delta-$perturbation $\La$ of $\Gamma$ the system  $E(\La)$ is
a Riesz basis for $L^2(\Ss)$.}

\medskip

This is proved in \cite{ya}, p. 161 for the case $\Ss=[-\pi,\pi]$. The proof of the general case is similar.

\subsection{Proof of Theorem 3.1}

          The property of exponential system to be a Riesz basis
         in $L^2(\Ss)$ is invariant with respect to translations
         of the set.  So, we may assume $x=0$ in the definition of $\Omega$.

         Denote by $\Omega(k), k=0,1,...,$ the following finite family of sets $\Ss$ in
         $\R^n$:
$$
                       \Ss=2^{-k}\bigcup_{j=1}^{2^{kn}}([0,2\pi]^n+2\pi m_j),
$$
         where the $m'$s are different vectors in $\Z^n, |m|< 2^k$.
         Observe that the sequence of these families is nested,
         and their union gives $\Omega$.

         Let $\epsilon>0$ be given. We construct $\La$ by the following induction
         process:          Set $\La_0 = \Z^n$ and $ \epsilon(0)=\epsilon/2.$
          Clearly, $E(\La_0)$ is a Riesz basis for  $\Ss\in \Omega(0)$.
          Suppose for $p=1,2,...,k-1$,  sets $\La(p)$ and positive numbers
          $\epsilon(p)$ are defined  such that:

            (i) $\La(p)$ is an $\epsilon(p) -$perturbation of $\La(p-1)$;

            (ii) $\La(p)$ is a union of $2^{np}$ shifted copies of $2^{p}\Z^{n}$;

            (iii) $E(\La(p))$ is a Riesz basis in $L^2(\Ss)$ for every $\Ss \in \Omega(p)$.

            (iv) $\epsilon(p)<\epsilon(p-1)/2 $.

             Describe the $k- $th  step of the induction.
          According to Lemma~3.2, we can fix
                    $\epsilon(k)< \epsilon(k-1)/2$
           such that every $\epsilon(k)-$perturbation $\La$ of $\La(k)$ is
           a Riesz basis in $L^2(\Ss)$ for every set $\Ss\in \Omega(k)$.

             Take $\Ss\in \Omega(k)$. Then Lemma 3.1
          (after re-scaling) provides a sufficient condition
          on vectors $u_1,u_2 ,... u_{2^{kn}}$
          in terms of the corresponding determinant, which ensures that
            $$\La(k):=\cup_{l=1}^{2^{kn}} (2^{k}\Z^{n} + u_l)$$ is such that $E(\La(k))$ is a Riesz basis
           in $L^2(\Ss)$. Then Lemma 3.3 and the classical uniqueness theorem
          for analytic functions imply that this condition is
          satisfied for almost every choice of vectors $\{u_l\}$.
          It follows that one can choose these
          vectors so that the set $\La(k)$ will satisfy (i)-(iv)
          with $p$ replaces by $k$.
           Finally, we get the set $\La$ as the pointwise limit of $\La(k)$ as $k\to\infty$.
           One can easily check  that it
          satisfies the requirements of the Theorem 3.1.

 \section{Sampling}


    In \cite{b3} Beurling obtained the following sufficient condition for sampling in  $B_\B$, where $\B$ denotes the unit ball in $\R^n$:

\medskip\noindent
{\bf  Theorem D}
    {\sl Assume  $\La\subset\R^n,n\geq1,$ satisfies
\begin{equation}\label{b1}
\La+\rho\B=\R^n \ \mbox{ for some } \rho<\frac{\pi}{2}.
\end{equation}  Then
 $\La$ is a sampling set for  $B_\B$. }

     \medskip

  One can check that condition (\ref{b1}) is equivalent to condition $\rho(\La)<\pi,$ where $\rho(\La)$ is defined in (\ref{rd}). So, for $n=1$  this theorem is equivalent to Theorem B stated in introduction.

   Theorem D is sharp for every $n\geq 1$ in the sense that  it ceases to be true for $\rho=\pi/2.$

   Beurling's approach in \cite{b3} is based on reduction of the
   sampling problem to possibility of
   balayage (or sweeping) of any given finite measure from $\R^n$ to the set $\La$ without changing the values on $\B$
   of its Fourier transform.

     Beurling noticed without proof that this approach
   works in a more general setting when the spectrum
   is a convex  body. See also Benedetto and Wu \cite{bw}, Theorem~7.2.

    Beurling  proved  a quantitative version of Theorem D:
    If (\ref{b1}) holds then
  $$
        \Vert f\Vert_\infty\leq \frac{1}{1-\sin \rho} \Vert f|_\La\Vert_\infty,\ \mbox{ for every } f\in B_{\B}.
$$
 He noticed without proof that the inequality above  can be replaced by a stronger one:
            \begin{equation}\label{bb}
     \Vert f\Vert_\infty\leq \frac{1}{\cos \rho}\Vert f|_\La\Vert_\infty, \ \mbox{ for every } f\in B_{\B} .
     \end{equation}

    We will present a completely different elementary
    approach, which provides this estimate in the
    general context of convex spectrum.

     Throughout this section  we assume that $\La$ is a set in $\R^n$ and that $\K$ is  closed convex central-symmetric body  of positive measure in $\R^n$. We also assume that the dimension $n>1.$  Recall that
     $\K^o$ denotes the dual body of $\K$.

\medskip\noindent
      {\bf  Theorem 4.1} {\sl Assume sets $\K$ and $\La$ satisfy \begin{equation}\label{lr}   \La+\rho\K^o=\R^n \ \mbox{for some } \rho<\frac{\pi}{2}.\end{equation}
            Then (\ref{bb}) is true,  and so
       $\La$ is an SS for $B_\K$. }

\medskip
Observe that Theorem D is a particular case of Theorem 4.1 in which $\K=\B$.

Theorems 4.1 and 2.1 imply

\medskip\noindent
      {\bf  Corollary 4.1} {\sl Assume $\La$ is a uniformly discrete  set. If
(\ref{lr}) is true then       $\La$ is an SS for $PW_\K$. }

\medskip

Theorem 4.1 ceasues to be true when $\rho=\pi/2.$ Moreover, sets $\La$ satisfying (\ref{lr}) with $\rho=\pi/2$ need not to be even uniqueness sets:

\medskip

\noindent {\bf Proposition 4.1} {\sl For every   $\K$ there exist $\La$  and a non-trivial function $f\in B_\K$ such that
$\La+\frac{\pi}{2}\K^o=\R^n$ and
 $f(\lambda)=0,\lambda\in\La$.}

   \subsection{Proof of Proposition 4.1}
Take points  $x_0\in \frac{\pi}{2}\K^o$ and $t_0\in \K$ such that $x_0\cdot t_0=\pi/2.$ The spectrum of the function $\sin(x\cdot t_0)$ consists of two points $\pm t_0\in\K,$ and so $\sin(x\cdot t_0)\in B_\K.$
Denote by $\La=\{x\in \R^n: x\cdot t_0\in \pi\Z\}$  the zero set of  $\sin(x\cdot t_0)$, and by
 $I=\{\tau x_0: -1\leq\tau\leq 1\}\subset \frac{\pi}{2}\K^o$ the interval from $-x_0$ to $x_0$. Clearly, for every point $y\in\R^n$ there exist $n\in\Z$ and $-1\leq\tau\leq 1$ such that
$y\cdot t_0 =\pi n-\tau \pi/2.$ Hence, $y-\tau x_0\in \La$, which implies $\La+I=\R^n$.

 \subsection{ Proof of Theorem 4.1}
We shall deduce Theorem 4.1  from the following

\medskip\noindent
{\bf Lemma 4.1} {\it Suppose a function $g\in B_{[-\tau,\tau]}$ satisfies $|g(0)|=\Vert g\Vert_\infty$. Then }
\begin{equation}\label{ch}|g(u)|\geq|g(0)|\cos(\tau u), \ \  |u|<\pi/2\tau.\end{equation}


This lemma is proved in \cite{cl} (see proof of Theorem 4). For completeness of presentation, we  sketch the proof below.

Let us prove Theorem 4.1. Take any function $f\in B_\K$. Assume first that $|f|$ attains maximum on $\R^n$, i.e. $|f(x_0)|=\Vert f\Vert_\infty$ for some $x_0\in\R^n$.  By (\ref{lr}),  there exists $\lambda_0\in\La$ with $\Vert\lambda_0-x_0\Vert_{\K^o}\leq\rho$.
  Consider the function of one variable $g(u):=f(x_0+u(\lambda_0-x_0)), u\in\R$.
 One may check that $g\in B_{[-\tau,\tau]}$ with $\tau=\Vert\lambda_0-x_0\Vert_{\K^o}$. Also, clearly $|g(0)|=\Vert g\Vert_\infty$ and $g(1)=f(\lambda_0)$. Since $\tau\leq\rho<\pi/2$, we may use inequality (\ref{ch}) with $u=1$:
 $$\Vert f\Vert_\infty=|f(x_0)|=|g(0)|\leq  \frac{|g(1)|}{\cos \tau}\leq
 \frac{|f(\lambda_0)|}{\cos \rho}\leq\frac{1}{\cos \rho}\Vert f|_\La\Vert_\infty. $$

If $|f|$ does not attain maximum on $\R^n$, we consider the function $f_\epsilon(x):=f(x)\varphi(\epsilon x)$, where
$\varphi\in B_{\epsilon\B}$ is any function satisfying $\varphi(0)=1$ and $\varphi(x)\to0$ as $|x|\to\infty.$
It is clear that $f_\epsilon\in B_{\K+\epsilon\B}$ and that $f_\epsilon$ attains maximum on $\R^n$. Set $g_\epsilon(u):=f_\epsilon(x_0+u(\lambda_0-x_0)), u\in\R,$ where $x_0$ and $\lambda_0$ are chosen so that $|g_\epsilon(0)|=\Vert f_\epsilon\Vert_\infty$ and  $\Vert\lambda_0-x_0\Vert_{\K^o}\leq\rho$. We have $g\in B_{[-\tau-\delta,\tau+\delta]} $,
 where $\tau=\Vert\lambda_0-x_0\Vert_{\K^o}\leq\rho<\pi/2$ and $\delta=\delta(\epsilon)\to0$ as $\epsilon\to0.$ So, if $\epsilon$ is so small that $\tau+\epsilon<\pi/2$, we may repeat the argument above to
obtain $\Vert f_\epsilon\Vert_\infty\leq \Vert f_\epsilon|_\La\Vert_\infty/\cos(\rho+\delta)$. By letting $\epsilon\to 0$, we obtain (\ref{bb}).

     \subsection{Proof of Lemma 4.1}

1. The proof in \cite{cl} is based on the following result from \cite{d} (for some extension see \cite{h}):
{\sl Let $f\in B_{[-\tau,\tau]}$ be a real function satisfying $-1\leq f(x)\leq 1$ for all $x\in\R$.
Then for every real $a$ the function $\cos (\tau z + a) - f(z)$ vanishes identically or else it
has only real zeros.  Moreover it has a zero in every interval where $\cos (\tau z + a)$ varies between -1 and 1 and all the zeros are
simple, except perhaps at points on the real axis where $f(x) = \pm  1.$}

Sketch of proof. We may assume $a=0$ and $\tau=1$. Consider the function
$$
 f_\epsilon(z):=(1-\epsilon)\frac{\sin (\epsilon z)}{\epsilon z}f((1-\epsilon)z).
$$
One may  check that  $f_\epsilon\in B_{[-1,1]}$, $-1<f(t)<1, t\in\R,$ and that the estimate holds
$$
|f_\epsilon (z)|\leq \frac{e^{ |y|}}{\epsilon |z|}, z=x+iy \in\C.
$$ This shows that $|f_\epsilon(z)|<|\cos z|$ when $z$ lies on a rectangular contour $\gamma$ consisting of segments of the lines $x = \pm N\pi, y= \pm N,$ where $N$ is every large enough integer. So by Rouch\`e's theorem the function $\cos z-f_\epsilon(z)$
has the same number of zeros in $\gamma$ as $\cos  z$, that is, $2N$ zeros. On the
real axis $|f_\epsilon|\leq 1-\epsilon $. Hence,  $\cos z-f_\epsilon(z)$ is alternately plus and
minus at the $2N+1$ points $k\pi$, $|k|\leq N,$ so  it has  $2N$ real zeros inside $\gamma$. Taking larger values of $N$ we see that $\cos z-f_\epsilon(z)$ has exclusively real and simple
zeros, which lie in the intervals $(k\pi,(k+1)\pi)$.

The zeros of $\cos z-f(z)$ are limit points of the zeros of $\cos z-f_\epsilon(z)$  as $\epsilon\to 0$. Thus $\cos  z-f(z)$ cannot
have non-real zeros. Moreover, it has an infinite number of real zeros
which are all simple, except those at the points $k\pi$ iff $f(k\pi) = ( - 1)^{k}.$
Every interval $k\pi<z<(k + 1)\pi$ at the endpoints of which
$| f (t) | < 1$ contains exactly one zero. If $f(k\pi) = (-1)^{k}$, we have a
double zero at $k\pi$ but no further zeros in the interior or at the endpoints
of the interval $((k - 1)\pi, (k+ 1)\pi). $

2. It suffices to prove Lemma 4.1 for functions $f\in B_{[-\tau,\tau]}$ that are real on $\R$. Since $f$ has a local maximum at $t=0$, the function $f(t)-\cos\tau t$ has a repeated zero at $t=0. $
By the discussion above we see that either $f(t)$ is identically equal to $\cos\tau t$ or $f(t)-\cos\tau t$ does not vanish on $[-\pi/\tau,0)\cup(0,\pi/\tau]$. It follows that  $f(t)>\cos\tau t$ on each of the intervals $[-\pi/\tau,0)$ and $(0, \pi/\tau]$.

   \section{Interpolation}

    In this section we discuss sharp sufficient condition
     for interpolation of discrete functions
     in $\R^n$ by functions with convex spectrum.
     We present here an extended version of our results briefly
     stated in  \cite{ou1}.

       \subsection{Ingham's theorem}

        The origin of the subject is Ingham's Theorem C from \cite{i} formulated  in introduction:

     \medskip\noindent{\bf Theorem C}  {\sl Let $\Ss$ be an interval and $\La\subset\R.$
        If  $\delta(\La) > \frac{2\pi}{|\Ss|}$  then $\La$ is an IS for $PW_\Ss$.}

        \medskip

     This result is sharp: It is proved in \cite{i} that the set
    $$\La=\{\pm \pi(n-\frac{1}{4}), n=1,2,3,...\}$$ is not an IS for $PW_{[-1,1]}$ (see also  Theorem 5, p. 103 in \cite{ya}),   so the condition $\delta(\La)\geq 2\pi/|\Ss|$
     is not sufficient for interpolation.

    We sketch the proof of Theorem C. The following  observation
    is important:
     Suppose there is a continuous even function $K(x)$, supported
     on   $\K:=(-r,r)$ such that   $   K(0)>0 , \hat K \in L^1(\R)$ and $\hat K$ is negative outside of $\Ss=[-1,1]$.
     Then $\delta(\La)>r$ implies that $\La$ is an IS for $PW_\Ss$.

         Indeed, denoting $b:= \max_{t\in\R}\hat K(t)$, we have for every polynomial   $P(t)=\sum_{\lambda\in\La}c(\lambda)e^{i\lambda t}$ that
$$
\int_{-1}^1|P(t)|^2\,dt\geq \frac{1}{b}\int_{-\infty}^\infty|P(t)|^2\hat K(t)\,dt=$$$$\frac{1}{b}\sum_{\lambda,\lambda'\in\La}c(\lambda)\bar c(\lambda')K(\lambda-\lambda')=\frac{K(0)}{b}\sum_{\lambda\in\La}|c(\lambda)|^2,
$$
    and due to (\ref{ip}) we get the result.

       Such a function $K$ does exist whenever $r>\pi$. Indeed, set
            $$
              K(x):= (1+D^2)(H\ast H), \ \ H(x):=\left\{ \begin{array}{ll} \cos (x\pi/r) &  |x|<r/2\\ 0 & |x|\geq r/2\end{array}\right.
  $$
            where $D=d/dx$ is the differentiation operator. One  can check that $K$ satisfies the properties above.

       Ingham-type results for $\R^n, n > 1$, were first obtained by  Kahane in  \cite{k1}, \cite{k}.

       Baiocchi, Komornik and Loretti \cite{bkl} extended Ingham's approach
     to the balls in $\R^n$. In this case one may take
$$
               K(x)= (1+\Delta)(H*H),$$      where $H$ is the first eigen-function of the Laplacian $\Delta$
      in $(r/2) \B$ with zero  boundary condition.
       This gives the following result (see \cite{bkl}, \cite{kl}):
      Every set $\La\subset\R^n$ satisfying $\delta(\La)>2\nu_n$ is an IS for $PW_\B$, where $\nu_n$ is the first root of the Bessel function $J_{n/2-1}$.

      Observe that $\nu_n=n(1+o(1)), n\to\infty.$ This shows that the linear growth $$\delta(\La)> 2n(1+o(1)), \ n\to\infty,$$ is sufficient for $\La$ to be an IS for $PW_\B$, where $\B\subset\R^n$. In the next section we show that the linear growth of $\delta(\La)$ is also necessary.


\subsection{Necessary condition}

          Let us consider the following general setting:  Let $\Ss$ and $\K$ be  two symmetric convex bodies in $\R^n$.
      The first body defines the Paley--Wiener space $PW_\Ss$, while the second one is used for measuring the distances between the interpolation nodes. To avoid trivial remarks, we assume that the dimension $n>1$.

          Consider the property:
\begin{equation}\label{is}
          (\La-\La) \cap \K =\{0\} \ \ \mbox{implies that $\La$ is an IS for $PW_\Ss$.}
\end{equation}
   Equivalently,
$$
             \La  \mbox{ is an IS for $PW_\Ss$ whenever}
             \Vert\lambda-\lambda^*\Vert_\K  > 1, \lambda,\lambda^*\in\La, \lambda\ne\lambda^*.
$$

\medskip\noindent{\bf Theorem 5.1} {\sl   If $\Ss$ and $\K$ satisfy (\ref{is}) then}
\begin{equation}\label{pn}
                  |\Ss|\cdot |\K| > (2\pi)^{n}.
\end{equation}

       When $\Ss=\B$ is the unite ball and  $\K=r\B$, we get

\medskip\noindent{\bf Corollary 5.1}
             {\sl If (\ref{is}) holds with $\Ss=\B$ and $\K=r\B$   then}
  $$
      r> \frac{2\pi}{|\B|^{2/n}}=\frac{n}{e} \, (1+o(1)),\ n\to\infty.
$$

          Hence, the linear growth of the "Ingham bound" for $\delta(\La)$
       is necessary for interpolation in $PW_\B$, while the "Beurling bound" does not depend on the dimension, see sec.4.

             \medskip
           The proof of Theorem 5.1 is based on
           special lattices described below.

       \subsection{Minkowski's lattices}

            The following problem goes back to Gauss.
         Let $\Q=[-1,1]^n$ be the unite cube in $\R^n$. Which part $w(n)$ of $\Q$ one
         can cover by disjoint balls of fixed radius?
         There is a huge number of papers devoted to the problem,
         and for many $n$ the optimal packing has been found.

          Notice that the simplest strategy, to decompose $\Q$ on small
         sub-cubes and to inscribe a ball in each one, gives the
         estimate  $w(n)\approx n^{-n/2}$.
         A much better estimate, which is asymptotically sharp,
         is due to Minkowski:
$$
                    w(n) > c^n ,
$$
          where $c$ is an absolute positive constant.
          For such a packing one should distribute the centra of balls
         in the vertexes of a special lattice.

           Minkowski proved that there is a lattice $\La$ in $\R^n$
           satisfying the following properties:
$$
                        \delta(\La) = 1, \ \ \mbox{ det}\,(\La)<|\B|,
  $$     where det$(\La)$ is the volume of it's fundamental parallelepiped, see  \cite{ca}, \cite{r}.
         An effective construction of such a lattice is not known
         so far. Minkowski's proof is based on ingenious  combination
         of algebraic and probabilistic argument.
            Later  Hlavka  extended his argument. Precisely,
         the following proposition is true, see  \cite{ca}, \cite{r}:
          Given a  convex symmetric body $\K$   there is a lattice $\La$ such that:
$$ \La\cap\K =\{0\}, \ \
                     \mbox{ det}(\La) < |\K|.$$

        \subsection{Proof of Theorem 5.1}

   To prove Theorem 5.1, take a lattice $\La$ satisfying $\La\cap\K =\{0\}$ and  det$(\La) < |\K|$. Since $\La-\La=\La$, we have $(\La-\La)\cap\K =\{0\}$, and so $\La$ is an IS for $PW_\Ss$. It follows from Proposition 3.3 (ii) that $\Ss+2\pi\La^*=\R^n$. Let $D$ be a fundamental parallelepiped of the lattice $2\pi\La^*$. Since $2\pi\lambda^*$-translations of $D$ tile $\R^n$, for almost every $x\in D$ there exists $\lambda^*\in\La^*$ with $x-2\pi\lambda^*\in \Ss$, and so we have
     $$|\Ss|\geq |D|=(2\pi)^n\, \mbox{det}\,(\La^*)=(2\pi)^n/\mbox{det}\,(\La)>(2\pi)^n/|\K|.$$

 \medskip\noindent{\bf Remark 5.1
 }  Theorems 5.1   admits extension to the  case when  $\Ss$ and $ \K$ are arbitrary (not necessarily convex) bounded  sets in $\R^n$.

 \medskip  The only change in the proof is that  uses a more general version
            of Minkovskii--Hlavka theorem for non-convex sets, see~\cite{ca}.

        \subsection{Sufficient condition}

          Here we prove a sufficient condition for interpolation
          sets. It is based on a concentration property for  functions
          and their Fourier transforms.
          We shall use it in the following form:

               \medskip\noindent{\bf
        Definition 5.2} {\sl  We say that a  couple     $(\Ss, \K)$
            admits concentration if there a function $F\in L^2(\R^n)$ satisfying $\Vert F\Vert_2=1$ and}
            \begin{equation}\label{cc}
\Vert F\Vert_{L^2(\K)}>\frac{7}{8}, \ \Vert \hat F\Vert_{L^2(\Ss)}>\frac{7}{8}.
\end{equation}

\medskip

Recall that in this paper we use a unitary form of the Fourier transform, so that $\Vert F\Vert_2=\Vert\hat F\Vert_2$.

         \medskip\noindent{\bf Theorem 5.2} {\sl    If a couple $(\Ss,\K)$ admits concentration
         then every set $\La$ satisfying $(\La-\La)\cap 12\K=\{0 \}$ is an IS for  $PW_\Ss.$}

                    \medskip\noindent{\bf Proof of Theorem 5.2}. 1. Assume that $(\Ss,\K)$ admits concentration.
            Take a function $F$ satisfying (\ref{cc})  and write $F=F_1+F_2$, where $F_1(x):=F(x)\cdot {\bf 1}_\K(x).$
            Then $\Vert \hat F\Vert_{L^2(\Ss)}> (7/8) \Vert  F\Vert_2$ and $\Vert F_2\Vert_2<(1/8)\Vert F\Vert_2,$ which gives
            $$
            \Vert \hat F_1\Vert_{L^2(\Ss)}\geq \Vert  \hat F\Vert_{L^2(\Ss)}-\Vert \hat F_2\Vert_2>\frac{3}{4}\Vert\hat F\Vert_2.
            $$
            Set
$ g(t):=\hat F_1(t)/\Vert \hat F_1\Vert_2.$ Clearly, we have
\begin{equation}\label{1}
                    g\in PW_\K, \ \Vert g\Vert_2=1, \ \Vert g\Vert_{L^2(\Ss)} >\frac{3}{4}.
\end{equation}

\medskip\noindent
2. Set
  $$A(t):=\left\{ \begin{array}{ll} 1 &  t\in 2\Ss\\ -1 & t\in \R^n\setminus 2\Ss \end{array}\right.$$  and consider the function
$$
h(t):=\int_{\R^n}|g(t-s)|^2 A(s)\,ds, \ \ t\in \R^n.
$$
               Clearly,
$$
                |h(t)| \leq 1, \  t\in\R^n.
$$                 Moreover, we have
\begin{equation}\label{4}
               h(t) > \frac{1}{2},  t\in \Ss, \
               h(t) <-\frac{1}{2},  t \in \R^n\setminus 3\Ss.
\end{equation}
              Indeed, let $\Vert t\Vert_\Ss\leq 1.$
                             Then, by (\ref{1}), $$h(t)=\int_{\R^n}A(t-s) |g(s)|^2\,ds =\int_{\Vert t-s\Vert_\Ss\leq 2} |g(s)|^2\,ds  - \int_{\Vert t-s\Vert_\Ss>2}|g(s)|^2\,ds\geq$$$$
                \int_\Ss|g(s)|^2\,ds -\int_{\R^n\setminus\Ss}|g(s)|^2\,ds >\frac{1}{2}, $$which gives the first inequality in (\ref{4}).
              The second one can be checked the same way.

 \medskip\noindent
3. Since $g\in PW_\K$ and $\Vert g\Vert_2=1$, it follows from Theorem A (i) that $g$ is an entire function satisfying (1) with $C=1$. Hence, the function $\bar g(\bar z)$ is entire and satisfies the same inequality.

 Observe now (see \cite{sw}, Lemma 4.11) that every function $\varphi\in PW_\K$ satisfies
$$\Vert\varphi(u+iv)\Vert_2:=(\int_{\R^n}|\varphi(u+iv)|^2\,du)^{1/2}\leq \Vert \varphi\Vert_2 e^{\Vert v\Vert_{\K^o}}.$$
Using this, it is easy to check that the formula$$
h(u+iv):=\int_{\R^n}g(u-s+iv)\bar g(u-s-iv) A(s)\,ds, \ \ u,v\in \R^n,
$$
defines an analytic continuation of $h$ into $\C^n$, and we have
  \begin{equation}\label{3}
|h(u+iv)|\leq \Vert g(u+iv)\Vert_2\Vert\bar g(u-iv)\Vert_2\leq e^{2\Vert v\Vert_{\K^o}}, \ \ u,v\in\R^n.
\end{equation}

\medskip\noindent 4.  Consider the function $ k(z): = h(z) g(z) \bar g(\bar z), z\in \C^n$.
              Clearly, $k$ is a bounded  and integrable on $\R^n$.   Also,  since $g\in PW_\K$, it follows from Theorem A (i) and (\ref{3}) that
$                        k \in PW_{4\K}.          $  So, the inverse  Fourier transform $K$ of $k$ is continuous and vanishes outside $4\K$.
                Due to (\ref{1}) and (\ref{4}) we have:
$  k(t) < 0, t  \in \R^n\setminus 3S,$ and
$$
                    K(0)=\int_{\R^n} k(t)\,dt \geq \frac{1}{2}\int_{\Ss}|g(t)|^2\,dt-\int_{\R^n\setminus \Ss}|g(t)|^2\,dt >0.$$
Repeating Ingham's argument involving the function $K$ (see sec. 5.1), we conclude  that every set $\La$ with $(\La-\La)\cap 4\K=\{0\}$ is an IS for $PW_{3\Ss}$. This  is equivalent to saying that every set $\La$ with $(\La-\La)\cap 12\K=\{0\}$ is an IS for $PW_{\Ss}$.

      \subsection{ Examples}

      Theorem 5.2 raises the following problem, which may have intrinsic interest: {\sl
      Which couples of convex bodies $\Ss,\K$ admit the concentration property?}

       Theorems 5.1 and 5.2   imply that condition $|\Ss||\K|>C^n$ is necessary, where $C>0$ is an absolute constant.       We shall give examples showing that in some cases      this condition is sharp.

       Denote by $\B_p\subset\R^n$  the unite ball in $l^p$-norm,
     that is
$$
        \B_p:= \{x=(x_1,...,x_n) , \Vert x\Vert_p^p:= \sum_{j=1}^n|x_j|^p\leq 1\}.
$$
  We shall keep the notation $\B$ for $\B_2$ and set $\Q:=\B_\infty$.

       We start with a simple

\medskip\noindent{\bf    Example 5.1} {\sl There is a constant $C$, which does not depend
        on $n$
        such that the couple $(\B,Cn\B)$ admits concentration.}

\medskip
            Indeed, consider the  Gaussian function $f(t)=\exp(-|t|^2/2)$.
       It is clear that given $\epsilon>0$ one can find $C>0$
       such that
$$
                \int_{|t|>Cn^{1/2}}  f^2(t)\,dt  <   \epsilon.
$$
        Since this function is invariant with respect to the
        Fourier transform, we get the concentration property for the couple
$              ( Cn^{1/2} \B ,Cn^{1/2} \B),$
        and the result follows.



      \medskip\noindent{\bf    Example 5.2} {\sl There is a constant $C$, which does not depend
        on $n$
        such that the couple $(\Q,Cn^{1/2}\B)$ admits concentration.}

\medskip

          This can be obtained from

                    \medskip\noindent
{\bf Lemma 5.1}. {\it Set
$$
h(x):=\prod_{j=1}^n\left(\frac{\sin  x_j}{ x_j}\right)^2, \ x=(x_1,...,x_n), \ n>1.
$$Then  $h\in PW_{2\Q}$, and  there exists $C>0$ such that }
$$
\int_{(C\sqrt{n/\epsilon})\B}|h(x)|^2\,dx >(1-\epsilon)\Vert h\Vert_2^2, \mbox{ for every } 0<\epsilon<1.
$$

\medskip \noindent{\bf Proof}. Set
$$
\beta:=\int_\R\left(\frac{\sin t}{t}\right)^4\,dt, \ \gamma:=\int_\R \frac{\sin^4 t}{t^2}\,dt.
$$
Then $ \Vert h\Vert_2^2=\beta^n$ and we have for every $R>0$ that
$$\int_{|t|>R}h^2(t)\,dt<\frac{1}{R^2}\int_{\R^n}|t|^2h^2(t)\,dt=\frac{n\gamma\beta^{n-1}}{R^2}.
$$
         Choosing $R= \sqrt{n \gamma/(\beta\epsilon)}$         we get the lemma.

   Clearly,  function $h$ is the Fourier transform
        of the convolution ${\bf 1}_Q\ast {\bf 1}_Q$ (up to a multiplicative
         constant), so that $h\in PW_{2\Q}.$
        This shows that        the couple $(\Q,Cn^{1/2} \B)$ admits
        concentration, for  some absolute constant $C$.

        By Lemma 5.1 and Theorem 5.2, there exists $C$ such that the  sets $\Ss=Cn^{1/2}\B$ and $\K=\Q$ satisfy
        condition (\ref{is}).      Using (\ref{pn}), one can easily check that this result
        is sharp in the sense that these sets do not satisfy (\ref{is}) when the constant $C$  is small. This also shows that if $C$     is  small enough, then the couple $(\Q,Cn^{1/2}\B)$ does not admit concentration.

        A more general form of Example 5.2 follows:

\medskip\noindent{\bf            Theorem 5.3} {\sl There is an absolute constant $C$
           such that the couple $             (n^{1/p} \B_p , C n^{1/q} \B_q)$
           admits concentration, provided  $p,q \geq1$ and at least one of the numbers does not exceed 2.}

\medskip\noindent{\bf           Proof. }
             Consider the family of convex bodies
$$
                         T_p := n^{1/p} \B_p.
$$
           Using Holder inequality one can prove the embedding
$                      T_p \supset T_q$, where $1\leq p\leq q\leq\infty.$
          Clearly, if $(\Ss,\K)$ admits concentration then every couple
          $(\Ss',\K')$ with $ \Ss\subset \Ss' ,\K\subset \K'$ also does.  Hence, Theorem 5.3 follows from Example 5.2.

\medskip\noindent{\bf Remark 5.2}  One can show that the Theorem 5.3 remains true
       for all  $p,q\geq 1$ provided one replaces $C$ by $C \min\{p,q\}$.

    \medskip   Since $\B_p^o=\B_q$, where $q=p/(p-1)$, Theorem 5.3 implies

\medskip\noindent{\bf   Corollary 5.2} {\sl
            There is an absolute constant $C$ such that  for every $p\geq 1$
       the couple $(\B_p,Cn\B_p^o)$ admits concentration.  }

\medskip\noindent{\bf             Remark 5.3}  The
           volumes of all bodies $T_p$ are comparable, up to
           the factor $C^n$. Hence, Theorem 5.1 shows that Theorem~5.3 is  sharp in the sense that it ceases to be true if the
           constant $C$ is small enough.

           \medskip
            In particular,
            in Corollary 5.2 we have  $C>1/e$, as it follows from
           Theorem 5.1 and the Santalo inequality that
$$
                        |\B_p||n\B_q| < (2\pi e +o(1))^n, \ n\to\infty.
$$

\medskip\noindent{\bf             Remark 5.4}  Using Theorem 2.1 we see that  Theorem 5.2 remains true for the Bernstein spaces $B_\Ss$.  The same is true for Theorem 5.1, provided one replaces condition (\ref{is}) with the condition
$$
          (\La-\La) \cap \K =\{0\} \ \ \mbox{implies that $\La$ is an IS for $B_\Ss$.}
$$

   \subsection{Questions}

    1.   The results above were first stated in our note \cite{ou1}.
         A general version of Corollary 5.2
       was also conjectured there:
      {\sl
        There is an absolute constant $C$ such that for every
       convex symmetric body $\Ss$ the couple
$ (\Ss, Cn\Ss^o) $
       admits  concentration, and therefore satisfies condition
       (\ref{is}).}

       We do not know whether the latter is true even for the
      lattices $\La$.  Due to Proposition 3.2, the question in
      this particular case can be  stated in purely  geometric form: {\sl
         Is it true that  if $\K+\Z^n$ is a "packing"
      then $Cn\K^o +\Z^n$ is a "covering" of $\R^n$?}

   2.  One can prove that for every constant $C$ the couple $(\Q,C\Q)$ does not
       admit concentration for all large $n.$    What about the condition (\ref{is})?

   3.  Is  condition (\ref{is}) symmetric,
       that is both pairs $(\Ss,\K)$ and $(\K,\Ss)$ satisfy it simultaneously?

\end{large}

\end{document}